\documentclass{article}%
\usepackage{amsmath}%
\usepackage{amsfonts}%
\usepackage{amssymb}%
\usepackage{graphicx}

\begin{document}

\title{On Relations Between the Stirling Numbers of First and Second Kind}
\author{Henrik Stenlund\thanks{The author is with Visilab Signal Technologies in Finland and is grateful for receiving support for this work.}
\\}
\date{March 25th, 2019}
\maketitle
\begin{abstract}
Four new relations have been found between the Stirling numbers of first and second kind. They are derived directly from recently published relations.\footnote{Visilab Report \#2019-04}
\subsection{Keywords}Stirling numbers of first and second kind
\subsection{Mathematical Classification}
Mathematics Subject Classification 2010: 11B73
\end{abstract}
\section{Introduction}
Only a few relations are known to exist between the Stirling numbers of first and second kind. Recently two new equations were introduced for this purpose. The Stirling numbers are extensively studied for instance in the field of combinatorics \cite{Jordan1933}, \cite{Bailey1935}, \cite{Boyadzhiev2009}. \cite{Noori2010} has displayed a nice historic review of the Stirling numbers.

The equations below seem to be the only nontrivial relations between the Stirling numbers of first and second kind \cite{Gradshteyn2007}.
\begin{equation}
S_{1\ n}^{(m)}=\sum_{k=0}^{n-m}{(-1)^{k}{\big(}^{n-1+k}_{n-m+k}{\big)}{\big(}^{2n-m}_{n-m-k}{\big)}S_{2\ {n-m+k}}^{(k)}} \label{eqn10}
\end{equation}
\begin{equation}
S_{2\ n}^{(m)}=\sum_{k=0}^{n-m}{(-1)^{k}{\big(}^{n-1+k}_{n-m+k}{\big)}{\big(}^{2n-m}_{n-m-k}{\big)}S_{1\ {n-m+k}}^{(k)}} \label{eqn15}
\end{equation}
\begin{equation}
{\delta}_{j,k}=\sum_{l=0}^{max(k,j)+1}{S_{1\ {l}}^{(j)}S_{2\ {k}}^{(l)}} \label{eqn20}
\end{equation}
\begin{equation}
{\delta}_{k,j}=\sum_{l=0}^{max(k,j)+1}{S_{1\ {k}}^{(l)}S_{2\ {l}}^{(j)}} \label{eqn25}
\end{equation}
The following two equations were recently found by the author \cite{Stenlund2019}
\begin{equation}
1=\sum_{j=1}^{m}{S_{1\ m}^{(j)}\sum_{k=1}^{j}{S_{2\ j}^{(k)}}} \label{eqn40}
\end{equation}
\begin{equation}
1=\sum_{j=1}^{m}{S_{2\ m}^{(j)}\sum_{k=1}^{j}{S_{1\ j}^{(k)}}} \label{eqn45}
\end{equation}
In the following, one is using mostly arithmetic operations in order to derive new expressions. Therefore it is not necessary to make too strict assumptions on number domains and differentiability etc. In general, the applied indexes $\in{N^+}$ and the variable $x\in{C}$ for instance.
\section{Stirling Numbers}
$S_{1\ i}^{(m)}$ is the Stirling number of first kind. The recurrence relation is the following
\begin{equation}
S_{1\ n+1}^{(m)}=S_{1\ n}^{(m-1)}-nS_{1\ n}^{(m)} \label{eqn120}
\end{equation}
and the special values are
\begin{equation}
S_{1\ n}^{(0)}=\delta_n^0, \ \ \ S_{1\ n}^{(n)}=1, \ \ \  S_{1\ n}^{(1)}=(-1)^{n-1}(n-1)!
\end{equation}
 $S_{2\ i}^{(m)}$ is the Stirling number of second kind. The recursion relation is
\begin{equation}
S_{2\ i+1}^{(m)}=mS_{2\ i}^{(m)}+S_{2\ i}^{(m-1)} \label{eqn155}
\end{equation}
and the special values are
\begin{equation}
S_{2\ n}^{(0)}=\delta_n^0,  \ \ \ S_{2\ n}^{(1)}=1,  \ \ \  S_{2\ n}^{(n)}=1 \label{eqn160}
\end{equation}
For both kinds of numbers, they are determined from the recursion relations above. Numbers outside the number triangles are zero. The basic features of the Stirling numbers of first and second kind can be found in most handbooks, see \cite{Abramowitz1970}, \cite{Gradshteyn2007}, \cite{Jeffrey2008}.
\section{Identities For the Stirling Numbers}
The equations below were found by the author \cite{Stenlund2019} while establishing the equations (\ref{eqn40}) and (\ref{eqn45}).
\begin{equation}
x^m=\sum_{j=1}^{m}{S_{1\ m}^{(j)}\sum_{k=1}^{j}{S_{2\ j}^{(k)}x^k}} \label{eqn1030}
\end{equation}
\begin{equation}
x^j=\sum_{m=1}^{j}{S_{2\ j}^{(m)}\sum_{k=1}^{m}{S_{1\ m}^{(k)}x^k}} \label{eqn1040}
\end{equation}
The parameter $x$ can have any value whatsoever. At the point $x=1$ one will get the results mentioned above. However, to derive new results one must start from the equations (\ref{eqn1030}) and (\ref{eqn1040}).
The summation (\ref{eqn1030}) can be opened up to the following polynomial
\begin{equation}
x^m=S_{1\ m}^{(1)}S_{2\ 1}^{(1)}x+S_{1\ m}^{(2)}[S_{2\ 2}^{(1)}x+S_{2\ 2}^{(2)}x^2]+
\end{equation}
\begin{equation}
...+S_{1\ m}^{(m)}[S_{2\ m}^{(1)}x+S_{2\ m}^{(2)}x^2+...+S_{2\ m}^{(m)}x^m]
\end{equation}
By using the known properties of the Stirling numbers one recognizes that the highest powers $x^m$ are eliminated from this equation leaving only powers lower than $m$.
\begin{equation}
0=\sum_{j=1}^{m-1}{S_{1\ m}^{(j)}\sum_{k=1}^{j}{S_{2\ j}^{(k)}x^k}}+S_{1\ m}^{(m)}\sum_{k=1}^{m-1}{S_{2\ m}^{(k)}x^k} \label{eqn1080}
\end{equation}
This leads at $x=1$ to
\begin{equation}
-\sum_{j=1}^{m-1}{S_{1\ m}^{(j)}\sum_{k=1}^{j}{S_{2\ j}^{(k)}}}=\sum_{k=1}^{m-1}{S_{2\ m}^{(k)}} \label{eqn1090}
\end{equation}
Identical steps starting from (\ref{eqn1040}) lead to
\begin{equation}
0=\sum_{m=1}^{j-1}{S_{2\ j}^{(m)}\sum_{k=1}^{m}{S_{1\ m}^{(k)}x^k}}+S_{2\ j}^{(j)}\sum_{k=1}^{j-1}{S_{1\ j}^{(k)}x^k} \label{eqn1180}
\end{equation}
Correspondingly, at $x=1$ one obtains
\begin{equation}
-\sum_{m=1}^{j-1}{S_{2\ j}^{(m)}\sum_{k=1}^{m}{S_{1\ m}^{(k)}}}=\sum_{k=1}^{j-1}{S_{1\ j}^{(k)}} \label{eqn1190}
\end{equation}
By differentiating equations (\ref{eqn1080}) and (\ref{eqn1180}), one will get at $x=0$ immediately
\begin{equation}
S_{2\ m}^{(1)}=-\sum_{j=1}^{m-1}{S_{1\ m}^{(j)}S_{2\ j}^{(1)}} \label{eqn1280}
\end{equation}
\begin{equation}
S_{1\ j}^{(1)}=-\sum_{m=1}^{j-1}{S_{2\ j}^{(m)}S_{1\ m}^{(1)}} \label{eqn1480}
\end{equation}
\section{Discussion}
The Stirling numbers of first and second kind form Pascal-type triangles with rather complicated closed-form equations for the numbers themselves. Not so many relations exist between these two kinds of Stirling numbers. 

The process started from the original equations (\ref{eqn1030}) and (\ref{eqn1040}) containing the variable $x$ by studying them in detail. The results were evaluated at $x=1$. Further two equations were produced by differentiating the equations (\ref{eqn1080}) and (\ref{eqn1180}) and evaluating them at $x=0$ which clears out all its powers. It is obvious that more complicated relations between the Stirling numbers of first and second kind can be generated by differentiating the equations (\ref{eqn1080}) and (\ref{eqn1180}) and evaluating at some value. These two polynomial equations seem to be in a key position while generating new equations.

The results indicate that the the two kinds of Stirling numbers are connected in a complicated way and there exist several complex relations. Equations (\ref{eqn1080}), (\ref{eqn1090}), (\ref{eqn1180}), (\ref{eqn1190}), (\ref{eqn1280}) and (\ref{eqn1480}) are believed to be new. 

\end{document}